\documentclass{article}
\usepackage{graphicx}
\usepackage{amsthm}
\usepackage{amssymb, amsmath}
\usepackage{color}
\usepackage{soul,xcolor}%To allow strike-out of text
\usepackage{url}

\begin{document}
\title{Riemann's zeta function and the broadband structure of pure harmonics}
\author{Artur Sowa\\
Department of Mathematics and Statistics, University of Saskatchewan\\
106 Wiggins Road,
Saskatoon, SK S7N 5E6,
Canada \\
sowa@math.usask.ca}
%\date{today}
\maketitle
\newtheorem{definition}{Definition}%[section]
\newtheorem{theorem}{Theorem}%[section]
\newtheorem{proposition}{Proposition}%[section]
\newtheorem{lemma}{Lemma}%[section]
\newtheorem{corollary}{Corollary}%[section]
\newtheorem{algorithm}{Algorithm}%[section]
\newtheorem{conjecture}{Conjecture}%[section]

\begin{center}

 Abstract

 \end{center}
Let $a\in (0,1)$ and let $F_s(a)$ be the periodized zeta function that is defined as $F_s(a) = \sum n^{-s} \exp (2\pi i na)$ for $\Re s >1$, and extended to the complex plane via analytic continuation. Let $s_n = \sigma_n + it_n, \, t_n >0 $, denote the sequence of nontrivial zeros of the Riemann zeta function in the upper halfplane ordered according to nondecreasing ordinates. We demonstrate that, assuming the Riemann Hypothesis, the Ces\`{a}ro means of  the sequence $F_{s_n} (a)$ converge to the first harmonic $\exp (2\pi i a)$ in the sense of periodic distributions. This reveals a natural broadband structure of the pure tone. The proof involves Fujii's refinement of the classical Landau theorem related to the uniform distribution modulo one of the nontrivial zeros of $\zeta$.
\vspace{.2cm}

  \noindent KEYWORDS: Riemann's zeta function, Fourier series, broadband
  \vspace{.2cm}

  \noindent AMS classification:  11K36, 42A99, 42C99, 11M35, 11M06

\section{Introduction} \label{Introduction}

The periodized zeta function is defined via the Dirichlet series:
\begin{equation} \label{periodized_def}
F_s(a) = F(s, a) = \sum\limits_{k=1}^{\infty} \frac{e^{2\pi i ka}}{k^s},\quad a \in (0, 1],
\end{equation}
where $s = \sigma + it$ is a complex variable. When convenient we will identify the segment $(0,1]$ with the unit circle $\mathbb{T}$ via $a\mapsto \exp (2\pi i a)$. Note that whenever $\sigma >1$ the series converges absolutely and uniformly in $a$, and hence $F_s$ is continuous on $\mathbb{T}$.
Moreover, as is well known, for any fixed $a\in (0,1]$ the function $s\mapsto F(s,a)$ can be extended into the half plane $\sigma \leq 1$ via analytic continuation; the Riemann zeta function is the particular case $\zeta(s) = F(s,1)$. In this way it is also seen that $F_s(a)$ is well defined for all $s$ and all $a\in (0,1)$\footnote{Note that $F_1(1)$ is not defined as it corresponds to the pole of $\zeta$.}. Let $s_n = \sigma_n + it_n, n\in\mathbb{N}$, be the sequence of the nontrivial zeros of the Riemann zeta function in the upper half plane, ordered so that $0< t_n\leq t_{n+1}$. It is known that $0<\sigma_n<1$ for all $n$.
Motivated by numerical evidence, see Fig. \ref{Fig}, we hypothesise that
\[
\frac{1}{N}\sum\limits_{n=1}^N F_{s_n}(a) \xrightarrow[N \rightarrow\infty]{} e^{2\pi i a} \quad \mbox{ a.e. in}\,\, (0,1].
\]
 The purpose of this article is to discuss the broader ramifications of this hypothesis and bring to light additional supporting arguments. In particular, we demonstrate convergence in the sense of distributions on $\mathbb{T}$, see Theorem \ref{FT_Atoms} in Section \ref{broadband}. Our proof relies on a calculation of the Fourier coefficients of $F_s$ for $s$ in the critical strip, see Theorem \ref{F_transform_thm} in Section \ref{calc_CS}. It also involves a careful estimate pertaining to the distribution of nontrivial zeros of the Riemann zeta function that was obtained by Fujii, \cite{Fujii2}, under the assumption of the Riemann Hypothesis (R.H.). In the closing sections we briefly discuss some new perspectives in harmonic analysis that the highlighted phenomenon appears to suggest, see Subsection \ref{encryption}. We also describe the numerical experiment that is critically important to the proposed hypothesis, see Subsection \ref{numerics}.

%Finally, the author believes that further exploration of the phenomena at hand will open new perspectives in harmonic analysis and its applications; for additional comments see the closing of Section \ref{broadband}.

%The following metaphorical interpretation of this phenomenon seems interesting. It is a well known fact that musical instruments produce sound waveforms that are characteristic of them. One can imagine an ensemble of musicians $1,2,3,\ldots, n, \ldots$, each playing an instrument whose characteristic waveform is, say, $\Re F_{s_n}$. Such an ensemble would collectively perform the sound of the pure tone when playing the lowest C. Now suppose we modify the instruments by removing the pure tone component from each of them. Such an orchestra, when in full blast, would perform pure silence.

\section{The Fourier coefficients of $F_s$ in the critical strip} \label{calc_CS}

Recall that for any $f\in L_p[0,1]$ with $p\geq 1$ its Fourier coefficients are defined via the Lebesgue integral to be:
\[
\hat{f} (k) = \int\limits_0^1 f(a) e^{-2\pi ika}\,da , \quad\quad k \in \mathbb{Z}.
\]
Since $( e^{2\pi i k \, .\,})_{k \in \mathbb{Z}}$ furnishes an orthonormal basis in $L_2[0,1]$ it follows directly from (\ref{periodized_def}) that $F_s\in L_2[0,1]$ whenever $\sigma >1/2$ and, in addition,
\begin{equation}\label{F_transform_trivial}
\hat{F}_s(k) = \left\{\begin{array}{lc}
                  k^{-s}, & k \geq 1 \\
                  0, & k \leq 0
                \end{array}\right.
\,\, \mbox{ provided } \sigma > 1/2.
\end{equation}
 Moreover, in this case the series in (\ref{periodized_def}) converges to $F_s$ in the $L_2$ norm.%\footnote{This is the classical Riesz-Fisher Theorem, \cite{R_E_Edwards}.}

Our first goal is to calculate the Fourier coefficients of $F_s$ when $0<\sigma<1$. Of course, a successful calculation should reproduce (\ref{F_transform_trivial}) when $1/2< \sigma <1$.
We begin by invoking the Hurwitz zeta function, $\zeta(s, a)$, which is related to $F_s(a)$ via an explicit formula discussed below, (\ref{H2oscil}). It is at first defined in the half-plane $\sigma > 1$ by
\begin{equation} \label{Hurwitz_def}
\zeta(s, a) = \frac{1}{a^s} + \sum\limits_{n=1}^{\infty} \frac{1}{(n+a)^s},\quad a \in (0, 1].
\end{equation}
Note that the series on the right converges absolutely and uniformly and therefore defines a continuous function of $a \in [0,1]$, while the first term on the right is merely continuous in $(0,1]$. The analytic continuation of $\zeta(s,a)$ to the entire $s$-plane is given by (see Theorem 12.3 in \cite{Apostol}):
\begin{equation} \label{Hurwitz_contour}
\zeta(s, a) = \Gamma(1-s)I(s,a) \quad (\mbox{for all } s )
\end{equation}
where $\Gamma$ is the gamma function and
\begin{equation}\label{def_I}
I(s,a) =\frac{1}{2\pi i}\int\limits_{\mathcal{C}} \frac{z^{s-1}e^{az}}{1-e^z}dz
\end{equation}
is a contour integral over $\mathcal{C}= \mathcal{C}_1 \cup \mathcal{C}_2 \cup \mathcal{C}_3$ with
\[
\mathcal{C}_1: z = re^{-i\pi}, c\leq r < \infty;  \quad\quad
\mathcal{C}_2: z = ce^{i\theta}, -\pi\leq \theta \leq \pi; \quad\quad
\mathcal{C}_3: z = re^{i\pi}, c\leq r < \infty .
\]
In other words, the contour $\mathcal{C}$ runs along the real axis from negative infinity to the point $-c$, where $c$ is arbitrary as long as $0< c< 2\pi$,  then circles counterclockwise around the origin and finally goes back to negative infinity along the same path. Thus, $z^s$ means $r^se^{-i\pi s}$ on $\mathcal{C}_1$ and $r^se^{i\pi s}$ on $\mathcal{C}_2$. It is well known that $s \mapsto I(s,a)$ is an entire function of $s$ (for any fixed $a\in (0,1]$). In consequence $\zeta(s,a)$ is analytic for all $s$ except for a simple pole at $s=1$ with residue $1$. The argument used in the proof of these facts, see \cite{Apostol}, shows also that the function is smooth with respect to $a$ in $(0,1]$. Indeed, it suffices to observe that the integral converges uniformly with respect to $a\in [\epsilon, 1]$.

Next, we find the Fourier coefficients of $a\mapsto I(s,a)$ when $\sigma <1$ (i.e. in the critical strip as well as in the left half-plane). Recall that $c < 2\pi$. First, assume $k\neq 0$.  In this case (\ref{def_I}) yields
\begin{align}
\hat{I}(s,k) &= \frac{1}{2\pi i}\int\limits_0^1 \int\limits_{\mathcal{C}} \frac{z^{s-1}e^{az}}{1-e^z}dz \, e^{-2\pi ika}da =
\frac{1}{2\pi i}\int\limits_{\mathcal{C}} \frac{z^{s-1}}{1-e^z} \, \int\limits_0^1 e^{az} e^{-2\pi ika}da\, dz  \label{one}  \\
&= \frac{1}{2\pi i}\int\limits_{\mathcal{C}} \frac{z^{s-1}}{1-e^z} \, \frac{e^z - 1}{z - 2\pi i k}\, dz = \frac{1}{2\pi i}\int\limits_{\mathcal{-C}} \frac{z^{s-1}}{z - 2\pi i k}\, dz \label{two}\\
& = \underset{z = 2\pi i k}{\mbox{Res}} \, \frac{z^{s-1}}{z - 2\pi i k} = (2\pi i k)^{s-1} \label{three}
\end{align}
Note that the contour with reversed orientation, i.e. $\mathcal{-C}$, encircles the point $z = 2\pi i k$ counterclockwise. The integral converges because $\sigma <1$. Since the contour leaves the branching point $z =0$ of $z \mapsto z^{s-1}$ on the outside the integral is evaluated the same way as for a meromorphic integrand. Note that $z = 2\pi i k$ is the only singularity of the integrand inside the contour, so that the contour integral over $\mathcal{C}$ may be replaced by a contour integral over $|z-2\pi i k | = \epsilon$, hence (\ref{three})\footnote{ The latter integral may also be evaluated directly via an application of the generalized binomial formula
\[
\left( 1+ \frac{\epsilon}{2\pi i k} e^{i\theta}\right)^{s-1} = 1+ (s-1)\frac{\epsilon}{2\pi i k}e^{i\theta} + \binom{s-1}{2}\left(\frac{\epsilon}{2\pi i k}\right)^2e^{i2\theta} + \ldots
\]}.
The interchange of the two integrals in (\ref{one}) is allowed by the Fubini Theorem, because the double integral converges absolutely for every $s$. Indeed, it suffices to show that one of the iterated integrals converges absolutely. Let $z = x+ iy$; we have
\[
\int\limits_{\mathcal{C}}\int\limits_0^1  \left|\frac{z^{s-1}e^{az}}{1-e^z} \, e^{-2\pi ika}\right| da dx\leq
\int\limits_{\mathcal{C}} \left|\frac{z^{s-1}}{1-e^z}\right| \frac{e^x -1}{x}dx,
\]
and an analogous inequality with $dx$ replaced by $dy$. Consider
\[
\int\limits_{\mathcal{C}} \left|\frac{z^{s-1}}{1-e^z}\right| \frac{e^x -1}{x}dz.
\]
The integrand is a smooth function on the compact circle $\mathcal{C}_2$. Indeed, $(e^x -1)/x$ is smooth at $x=0$ while $|1-e^z|$ is bounded below by a positive constant. This means that this part's contribution is finite. On $\mathcal{C}_1$ and $\mathcal{C}_3$ we have $z=x=-r$ with $r\geq c$, and the integrand is bounded by $e^{\pi |t|}r^{\sigma -2}$ and, since $\sigma <1$, these integrals converge. This justifies the interchange of the integrals. In addition, the same argument shows that $a \mapsto I(s,a)$ is in $L_1[0,1]$; see also (\ref{singularity}) and the remark that follows.

 Second, consider the case $k=0$. Proceeding as in lines (\ref{one}) and (\ref{two}) one can represent $\hat{I}(s,0)$ as a contour integral over $\mathcal{C}$, then notice that the integrand is regular outside the contour, which shows that $\hat{I}(s,0) =0$. Alternatively, this can be observed via an explicit integration over the three parts of the contour.

By virtue of (\ref{Hurwitz_contour}) we also have $\zeta(s,.)\in L_1[0,1]$, and
\begin{equation}\label{H_transform}
 \hat{\zeta}(s,0) = 0,\,\, \hat{\zeta}(s,k) =\Gamma(1-s) (2\pi i k)^{s-1} \mbox{ for } k\in \mathbb{Z}\setminus\{0\}\quad\quad (\sigma <1).
\end{equation}
For clarity we emphasize the interpretation of the complex powers, namely:
\[
(2\pi i k)^{s-1} = \left\{\begin{array}{ll}
                            (2\pi k)^{s-1} \, e^{i\, \pi (s-1)/2}, & k >0 \\
                            (2\pi |k|)^{s-1} \, e^{-i\, \pi (s-1)/2}, & k<0.
                           \end{array}
 \right.
 \]
The function $a\mapsto \zeta(s,a)$ has a singularity at the left end-point of the unit interval. This singulary is described by the following estimate, see Theorem 12.23, \cite{Apostol}. First, on the right from the critical strip (i.e. for $\sigma>1$) the distance $|\zeta(s, a) - a^{-s}|$ is bounded by the constant $\zeta(\sigma)$. Second, if $1-\delta \leq \sigma \leq 2$ for some $\delta \in (0,1)$, then
\begin{equation}\label{singularity}
|\zeta(s,a) - a^{-s}| \leq A(\delta) \,|t|^{\delta}\quad \mbox{ for } |t|> 1,
\end{equation}
where the constant $A(\delta )$ depends on $\delta$ but not on $s$. This implies that in fact $\zeta(s,.)\in L_p[0,1]$ whenever $1\leq p < 1/\sigma$ with $0<\sigma <1$ (and $|t| >1$).

We are now in a position to deduce the following:

\begin{theorem} \label{F_transform_thm}   For an arbitrary $s\in\{ 0< \sigma<1,\, |t|>1\}$ we have
$ F_s \in  C^\infty(0,1)$ and $ F_s\in L_p[0,1]$ whenever $1\leq p < 1/(1-\sigma)$. Moreover,
\begin{equation}\label{F_transform}
\hat{F}_s(k) = \left\{\begin{array}{ll}
                         k^{-s},  & k\geq 1 \\
                        0, & k \leq 0
                      \end{array}\right.
                      \,\, \mbox{ provided }  \sigma >0.
\end{equation}
Finally,
\begin{equation}\label{F_ae}
F_s(a) = \sum\limits_{k=1}^\infty k^{-s} e^{2\pi i k a} \quad \mbox{ for all } a\in (0,1).
\end{equation}
\end{theorem}

\noindent \emph{Proof.} We invoke the Hurwitz formula\footnote{This formula is a direct consequence of Theorem 12.6 in \cite{Apostol}. It is given explicitly in \cite{NIST} and also in \cite{Knopp} albeit in the latter it is printed with incorrect sign.}:
\begin{equation}\label{H2oscil}
  F_s(a) = i (2\pi)^{s-1} e^{-\pi i s/2} \Gamma(1-s) \left\{\zeta(1-s,a) - e^{\pi i s} \zeta(1-s,1-a)\right\}\quad (s \in \mathbb{C},\, 0< a<1).
\end{equation}
  In light of (\ref{singularity}) it is now clear that $ F_s\in L_p[0,1]$ whenever $1\leq p < 1/(1-\sigma)$ with $0<\sigma <1$ (and $|t|>1$).

  In light of (\ref{F_transform_trivial}) it suffices to calculate the Fourier coefficients for $0<\sigma< 1$. These are obtained by substituting (\ref{H_transform}) into (\ref{H2oscil}), observing that if $g(a) = f(1-a)$ then $\hat{g}(k) = \hat{f}(-k)$, as well as making use of Euler's formula
$\Gamma(s) \Gamma(1-s) = \pi/\sin{(\pi s)}$. This proves (\ref{F_transform}).

 To prove (\ref{F_ae}) we make use of the fact that $F_s$ is smooth in $[\epsilon, 1-\epsilon]$ for small $\epsilon >0$. Let us fix an arbitrary $a\in (0,1)$. There is a neighbourhood $a\in [\epsilon, 1-\epsilon]$ where the function is of bounded variation. It is well known that this implies convergence of the Fourier series at $a$, see \cite{R_E_Edwards}, Chapter 10\footnote{Note that convergence of the Fourier series to the value of the function \emph{almost everywhere} in $(0,1)$ follows also from $L_p$ ( $p>1$ ) integrability by virtue of the celebrated theorems of L. Carleson and R. A. Hunt.}. $\Box$
\vspace{.5cm}

\noindent
\textbf{Remark 1.}  We bring to the reader's attention consistency of (\ref{F_transform}) with (\ref{F_transform_trivial}).
In a way we have come a full circle starting from definition (\ref{periodized_def}) and ending with an identically looking (\ref{F_ae}). While  convergence in (\ref{F_ae}) is  only conditional for $s$ is in the critical strip, it is a quite gratifying to see that the series returns the value of the function $F_s(a)$ everywhere in $(0,1)$. We note that $F_s$ is not only smooth in $(0,1)$ but, in fact, it is analytic. Indeed, consider
\[
\tilde{F}_s(z) = \sum\limits_{k=1}^\infty k^{-s} z^{k} \quad \mbox{ for complex } z.
\]
It is straightforward that for $s$ in the critical strip the radius of convergence of this series is $1$, which means that the function $\tilde{F}_s$ has at least one pole on the unit circle. On the other hand, by Abel's theorem $\tilde{F}_s(z) \rightarrow F_s(a)$ when $z\rightarrow \exp(2\pi i a)$ within a Stolz angle, provided $F_s(a)$ is finite. Thus, the only pole of $\tilde{F}_s$ occurs at $z=1$, and the function is analytic in the neighbourhood of any other point on the unit circle, where it coincides with $F_s$.
\vspace{.5cm}

\noindent
\textbf{Remark 2.} Recall that the convolution of two functions $f,g\in L_1[0,1]$ is defined via
 \[
 f \star g(x) = \int\limits_0^1 f(x-y) g(y) dy
 \]
 It is a classical fact, \cite{R_E_Edwards}, that $f\star g \in L_1[0,1]$ and $(f\star g)\,\hat{} \, (k) = \hat{f}(k) \, \hat{g}(k)$.
 Theorem \ref{F_transform_thm} together with (\ref{F_transform_trivial}) imply that the family of  functions $\{F_s: \sigma > 0 \}$ forms a \emph{semigroup} with respect to convolution, i.e.
 \[
 F_{s} \star F_{s'} = F_{s+s'} \quad \mbox{ whenever } \sigma, \sigma' >0.
 \]
\vspace{.5cm}

\noindent
\textbf{Remark 3.}
 We emphasize that the theorem does not resolve the case $\sigma < 0$, for which formula (\ref{F_transform}) would be manifestly false. Indeed, it may be seen via an estimate similar to (\ref{singularity}), see Theorem 12.23, \cite{Apostol}, and the Hurwitz formula (\ref{H2oscil}) that $F_s\in L_1[0,1]$ even for negative $\sigma $ (assuming, as usual, that $|t|>1$). However, the Fourier coefficients of an integrable function converge to zero as $|k| \rightarrow \infty$ (by the Riemann-Lebesgue lemma), contrary to what formula (\ref{F_transform}) would imply in this case. %This limitation of our method stems from the fact that (\ref{H_transform}) holds only for $\sigma < 1$.

\section{A broadband representation of the pure tone} \label{broadband}

We will apply some classical results concerning the distribution of the nontrivial zeros of the Riemann zeta function. Let $\{x\}$ denote the fractional part of $x$. Recall that a sequence $(x_n)_{n=1}^{\infty}$ is said to be uniformly distributed modulo one
if
\[
\lim\limits_{N\rightarrow \infty}\frac{1}{N} \# \{1\leq n \leq N: \{x_n\} \in [\alpha, \beta ) \} = \beta - \alpha \quad \mbox{ for all } [\alpha, \beta )\subset [0,1 ).
\]

Let $s_n = \sigma_n + it_n$ be the sequence of zeros of Riemann's zeta function in the upper half plane, arranged in the order of growing imaginary parts\footnote{Here, the zeros are listed with multiplicities. It is not known at present if all the nontrivial zeros of zeta are simple; see \cite{Montgomery} for a discussion of this problem.}, so that $t_n\leq t_{n+1}$. It has been known for quite some time now\footnote{The discovery of this theorem has an interesting history; for brief historical remarks see \cite{review}, and for a more extended historical commentary see \cite{Steuding}.},  \cite{Elliott}, \cite{Hlawka}, \cite{Rademacher}, that
for every real number $\alpha \neq 0$ the sequence $(\alpha \, t_n)_{n=1}^{\infty}$ is uniformly distributed modulo one. Closely related to this fact is the Landau theorem, \cite{Landau}, which states the following: For a real number $x>1$ we have
\[
\sum\limits_{0< t_n\leq T} x^{-s_n} =- \frac{T}{2\pi} \frac{\Lambda(x)}{x} + O(\log T)\quad \mbox{ as } T\rightarrow \infty
\]
where $\Lambda $ is the von 'Mangoldt function, i.e. $\Lambda (x) = \log p$ when $x = p^m$ for a prime $p$ and $m\in \mathbb{N}$, and $\Lambda(x) = 0$ otherwise. Recall also that the number of zeros of $\zeta$ with $t_n\in[0,T]$ is known to be, see e.g. \cite{H_M_Edwards},
\begin{equation}\label{NofT}
N(T) = \frac{T}{2\pi} \log\frac{T}{2\pi e} + O(\log T )\quad \mbox{ as } T\rightarrow \infty.
\end{equation}
This shows that
\[
\frac{1}{N(T)}\sum\limits_{0< t_n\leq T} x^{-s_n} \xrightarrow[T \rightarrow\infty]{} 0.
\]
However, this result gives no information about the dependence of the rate of convergence on $x$. Indeed, the error term in Landau's asymptotic formula depends on $x$. Fortunately, more recently the formula was refined by other authors who analyze how the error depends on $x$. Specifically,  we will apply a theorem of Fujii, \cite{Fujii2}, which is obtained with the assumption that the R.H. holds, i.e. $\sigma_n = 1/2$ for all $n$. For its exact (rather lengthy) formulation the reader is referred to the source, i.e. Theorem 2 in \cite{Fujii2}. Here, we summarize as follows:
\begin{equation}\label{means_control}
\mbox{ For any } \sigma >0,\quad \frac{1}{N(T)}\sum\limits_{0< t_n\leq T} x^{-\sigma - i t_n} = O\left(\frac{x^{1/2-\sigma}\log x}{\log T}\right) \quad \mbox{ as } T\rightarrow \infty,\, x\rightarrow \infty .
\end{equation}
%The advantage is in control of the rate of convergence (to zero as $T\rightarrow\infty$) on $x$.
With this understood we turn to a discussion of the main problem: the broadband decomposition of the pure harmonic.
%Note first that under the R.H. none of the functions $F_{s_n}$ are square integrable; indeed this is a direct consequence of (\ref{F_transform}) and the Parseval identity. At the same time we have the following:

\begin{theorem}\label{FT_Atoms}
Assume $\sigma > 0$. Under the R.H.
\begin{equation}\label{weak_convergence}
  \lim\limits_{T \rightarrow \infty} \frac{1}{N(T)}\sum\limits_{t_n\leq T} F_{\sigma + i t_n} = \exp{(2\pi i \, .\,)},
\end{equation}
where the convergence type is as follows:
\begin{enumerate}
\item
distributional whenever $\sigma > 0$,

\item
in $L_2$-norm whenever $\sigma > 1$,

\item
uniform whenever $\sigma > 3/2$.
\end{enumerate}
\end{theorem}

\noindent \emph{Proof.} We note that Theorem \ref{F_transform_thm} applies to each $F_{s_n}$ as indeed $t_n\geq t_1 \simeq 14.1347 >1$. Thus, $F_{\sigma + i t_n}\in L_{1}$. Moreover, in light of (\ref{F_transform}) and (\ref{means_control}) for sufficiently large $T$ we have
\begin{equation}\label{convergence}
 \frac{1}{N(T)}\sum\limits_{t_n\leq T} \hat{F}_{\sigma + it_n}(k) =  O\left(\frac{k^{1/2-\sigma}\log k}{\log T}\right)\quad \mbox{  whenever }    k > 1.
\end{equation}
To prove pt. \emph{1.} we recall that whenever $f\in L_1[0,1]$ and $g$ is of bounded variation in $[0,1]$, the Parseval formula holds, see \cite{R_E_Edwards}, Chapter 10, i.e.
\[
\langle f, g\, \rangle : = \int\limits_0^1 f(a) \, \overline{g(a)} da = \sum\limits_k \hat{f}(k) \,\overline{\hat{g}(k)} .
\]
Next, let us fix an arbitrary function $g\in C^\infty (\mathbb{T})$. The fact that $\hat{g}(k) = o(1/k^m)$ for arbitrarily large $m$ together with (\ref{convergence}) imply
\begin{equation}\label{main_argument_convergence}
\left\langle \frac{1}{N(T)}\sum\limits_{t_n\leq T} F_{\sigma + it_n} - \exp{(2\pi i \, .\,)} ,\, g\,\right\rangle = \sum\limits_{k>1} \frac{1}{N(T)} \sum\limits_{t_n\leq T} \hat{F}_{\sigma  + it_n}(k)  \,\,\overline{\hat{g}(k)}
\xrightarrow[T \rightarrow\infty]{} 0.
\end{equation}
This proves the first case.
To prove \emph{2.} it is enough to observe that in light of (\ref{convergence}) $\sigma > 1$ implies
\[
\left\|\frac{1}{N(T)}\sum\limits_{t_n\leq T} F_{\sigma + it_n} - \exp{(2\pi i \, .\,)} \right\|_2 = O\left(\frac{1}{\log T}\right)
\]
 Similarly, pt \emph{3.} follows from an observation that if $\sigma > 3/2$, then
 \[
\left\|  \frac{1}{N(T)}\sum\limits_{t_n\leq T} F_{\sigma + it_n} - \exp{(2\pi i \, .\,)} \right\|_\infty = O\left(\frac{1}{\log T}\right).
\]
This completes the proof. $\Box$
\vspace{.5cm}

\noindent
\textbf{Remarks on the proposed hypothesis.}
As mentioned in Section \ref{Introduction} it may be hoped that convergence in (\ref{weak_convergence}) is in reality stronger than Theorem \ref{FT_Atoms} indicates. Indeed, one may interpret Fig. \ref{Fig} as supporting the claim that \emph{sequence $(F_{s_n} )$ is Ces\`{a}ro convergent to the fundamental tone pointwise almost everywhere}. Note that convergence everywhere in $(0,1)$ cannot hold, e.g. there is no convergence at $a=1/2$  because $F_{s_n}(1/2) = 0$ (since the zeros of the alternating zeta function coincide with the zeros of $\zeta$ in the critical strip)\footnote{It has been brought to the author's attention that convergence is also known to fail at the rational points $a=k/l$ where $l$ is a square-free integer.}. Also, convergence in the $L_\infty$-norm seems unlikely. Another subtle question is whether or not the optimal convergence result is inextricably dependent on the R.H.
%\vspace{.5cm}

\section{Addenda}
%\noindent
\subsection{Broadband encryptions of periodic functions}\label{encryption}

Suppose $f$ is a sufficiently regular function on $(0, 2\pi]\equiv\mathbb{T}$. In fact, in order to focus attention we assume  $\sum_{m\in \mathbb{Z}}|\hat{f}(m)| < \infty$, so that the Fourier series of $f$ converges uniformly and therefore $f$ is continuous on $\mathbb{T}$. In order to simplify notation we also assume $\hat{f}(0) = 0$. In light of Theorem \ref{FT_Atoms} it is natural to consider the following \emph{encryption} of $f$ (depending on fixed $\sigma$ and $n$):
 \[
 M_{\sigma, n} [f](a) =  \sum\limits_{m=1}^\infty \hat{f}(m) F_{\sigma + it_n}(ma) +  \sum\limits_{m=1}^\infty \hat{f}(-m) \overline{F_{\sigma + it_n}(ma)}.
 \]
 Theorem \ref{F_transform_thm} shows that the frequency modes of $F_{\sigma + it_n}$ decay in magnitude as $k^{-\sigma}$, i.e. very slowly for $\sigma\in(0,1]$. Hence, $M_{\sigma, n}[f]$ is typically going to have a lot of high-frequency content, constituting a \emph{broadband} encryption of the function $f$.  Nevertheless, in light of Theorem \ref{FT_Atoms}, one can expect that much information about $f$ can be retained in a family of broadband models, and recovered when desirable via simple averaging. Indeed, one may expect
 \begin{equation}\label{heurist}
 \frac{1}{N(T)} \sum\limits_{n\leq N(T)} M_{\sigma, n}[f] \approx f.
 \end{equation}
 Let us quantify this statement, assuming as usual that $\sigma > 0$.  First, let us fix $\sigma$ and denote
\[
\Phi_T =  \frac{1}{N(T)}\sum\limits_{t_n\leq T} F_{\sigma + it_n}.
\]
Now, let $\Phi_{T,m}$ be defined by $\Phi_{T,m} (a) = \Phi_T (ma)$ for $m \in\mathbb{N}$, and also $\Phi_{T,-m} = \overline{\Phi_{T,m}}$. We have
 \[
 \frac{1}{N(T)} \sum\limits_{n\leq N(T)} M_{\sigma, n}[f]  - f = \sum\limits_{m\in\mathbb{Z}} \hat{f}(m) (\Phi_{T,m}(a) - e^{2\pi i ma}).
 \]
In essence the argument in (\ref{main_argument_convergence}) relies upon the fact that, for a smooth $g$, $\hat{\Phi}_T (k)\,\overline{\hat{g}(k)} = O(\log k /(k^{1+\varepsilon} \log T))$ whenever $k>1$.  Now, observe that $\hat{\Phi}_{T,m} (mk) = \hat{\Phi}_T (k)$ while the Fourier coefficients corresponding to indices not divisible by $m$ all vanish.  This implies $\hat{\Phi}_{T,m} (k)\,\overline{\hat{g}(k)} = O(\log k /(|m|^{1+\varepsilon}k^{1+\varepsilon} \log T))$ for $k>1$. Hence
\[
\left\langle  \frac{1}{N(T)} \sum\limits_{n\leq N(T)} M_{\sigma, n}[f]  - f ,\, g\,\right\rangle = \sum\limits_{m}\hat{f}(m)\sum\limits_{k>1}  \hat{\Phi}_{T,m} (k) \,\,\overline{\hat{g}(k)}
\xrightarrow[T \rightarrow\infty]{} 0,
\]
i.e. the double series on the right converges absolutely for sufficiently large $T$ and its limit is zero as $T$ approaches infinity.  This endows (\ref{heurist}) with one rigorous interpretation. Again, it is an open problem to describe the optimal type of convergence that takes place in this scenario. However, as in Theorem \ref{FT_Atoms}, it is clear that convergence becomes increasingly stronger as $\sigma$ increases. It is interesting to observe the inherent counterbalance of this phenomenon with the diminishing high frequency content in the encryptions as $\sigma$ increases.

In addition, it is interesting to observe that any attempt at recovery of signal $f$ from $M_{\sigma, n}[f]$ (i.e. \emph{decryption}) would be a complicated matter when the waveform $F_{\sigma + it_n}$ is not known a priori. However, when it is known and in addition its Fourier coefficients diminish relatively rapidly, then decryption is a simple matter, see \cite{sowa21}.  Finally, we point out that when analogous operations are considered for digital signals the encryption as well as the decryption can be computed with high numerical efficiency, \cite{sowa22}, although the problem of numerical stability of these processes would warrant a separate discussion in the case of slowly diminishing Fourier coefficients (small $\sigma$).

\vspace{.5cm}

\noindent
\subsection{Details of the numerical experiment}\label{numerics}

\begin{figure}[ht!]
\includegraphics[width=175mm]{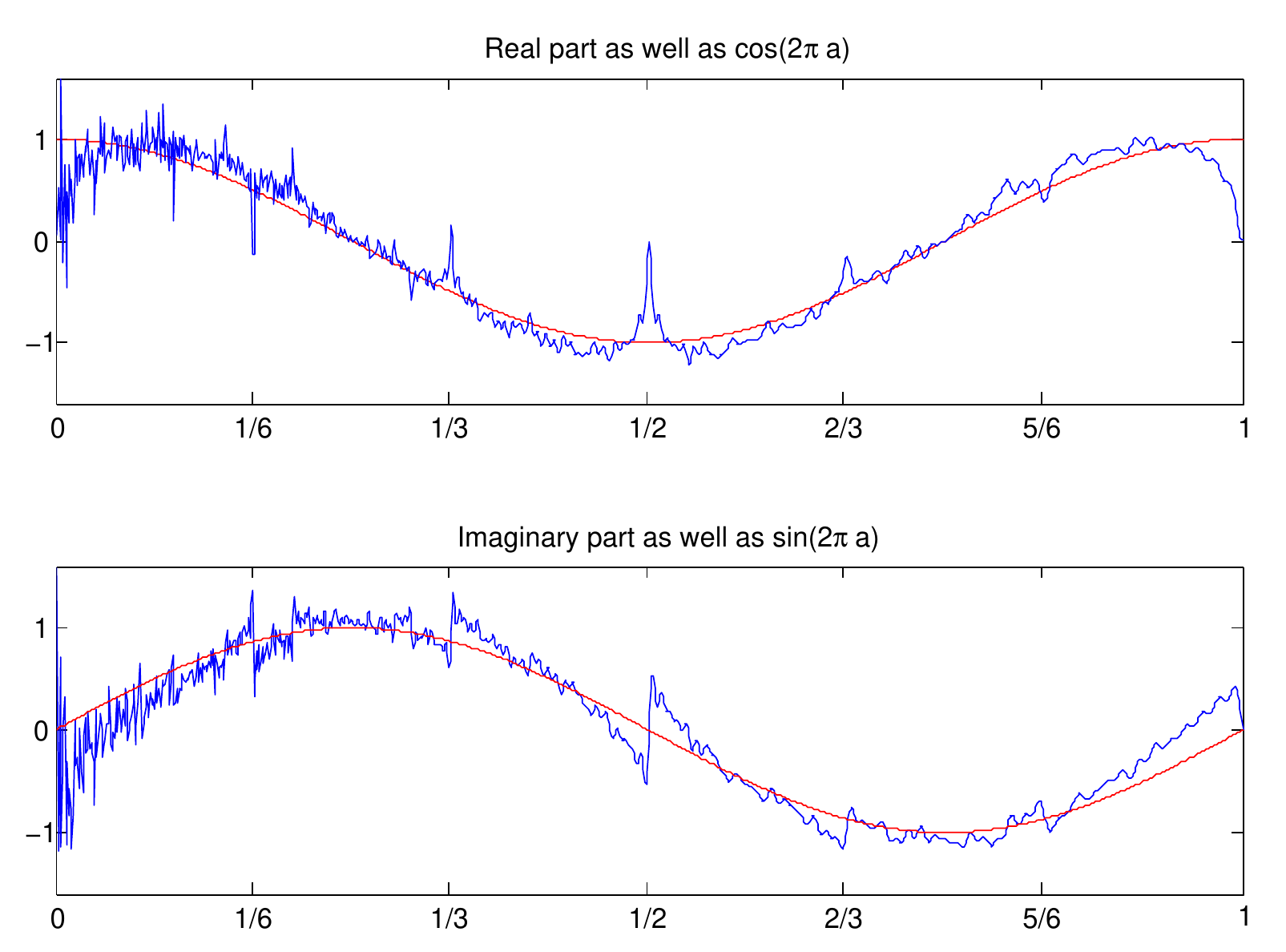}
%\emph{Fig.}
\caption{The real and imaginary parts of $\sum_{n\leq 237} F_{s_n}/237$ in comparison with the cosine and the sine functions. Note that the oscillatory patterns intensify toward the left end of the interval. This is related to the fact that $F_{s_n}$ break into chirpy oscillations at this end. Note also that the function assumes value zero in the middle of the interval. This is rigorously true because $s \mapsto F_{s}(1/2)$ coincides with the alternating zeta function and hence $F_{s_n}(1/2)=0$ for all $n$. }
\label{Fig}
\end{figure}

\vspace{1cm}
The graphs displayed in Fig. \ref{Fig} have been obtained by means of numerical evaluation of $F_s(a)$ for $\sigma >0$. It is based on the Hurwitz formula (\ref{H2oscil}) as well as the following classical estimate for $\zeta(s,a)$, see Theorem 12.21, \cite{Apostol}:
\[
\zeta(s,a) = \sum\limits_{n=0}^{N} \frac{1}{(n+a)^s} + \frac{(N+a)^{1-s}}{s-1} + o(N).
\]
Thus, the approximate values of $\zeta(s,a)$ are obtained from the first two terms on the right.
The rate of convergence of the $o(N)$ term depends on the size of $|s|$. Generally, one has to engage more steps (larger $N$) when $|s|$ is large,  especially so for $s$ in the critical strip, i.e. where the objects at hand are the most interesting. Also, for relatively large $|s|$ the computation of $\Gamma(1-s)$, an ingredient in (\ref{H2oscil}), becomes more problematic. For moderate magnitudes of $s$ qualitatively tenable examples can still be computed via asymptotic formulas.

For the numerical values of the sequence $(t_n)$ we have relied on data from  \cite{Odlyzko}. It is worthwhile pointing out that when $s = s_n$ the size of factor $e^{i\pi s}$ in (\ref{H2oscil}) is very small ($ \sim 5.1 \times 10 ^{-20}$ already for $s_1$). This implies that the difference between $F_{s_n}(a)$ and the function
\[
a \mapsto i (2\pi)^{s-1} e^{-\pi i s/2} \Gamma(1-s) \zeta(1-s,a), \quad s = s_n
\]
is so minuscule as to be negligible in most numerical applications, e.g. even for a fairly fine discretization of variable $a\in(0,1)$ the graphs $a\mapsto F_{s_n}(a)$ appear smooth at the right end-point (in spite of having a singularity  there of type $i e^{-\pi t_n} (1-a) ^{-.5 - i t_n}$).

%\vspace{1.5cm}

%%%%%%%%%%%%%%%%%%%%%%%%%%%%%%%%


\clearpage
\begin{thebibliography}{[90]}


\bibitem{Apostol} T. M. Apostol, Introduction to Analytic Number Theory, Springer-Verlag (1976)

\bibitem{NIST} T. M. Apostol, "Hurwitz zeta function", in F. W. J. Olver; D. M. Lozier; R. F. Boisvert; C. W. Clark, NIST Handbook of Mathematical Functions, Cambridge University Press (2010)

\bibitem{R_E_Edwards} R. E. Edwards, Fourier Series, A Modern Introduction, Vols. I and II,  Springer-Verlag (1979), second edition

\bibitem{H_M_Edwards} H. M. Edwards, Riemann's Zeta Function, Academic Press (1974)

  \bibitem{Elliott} P. D. T. A. Elliott, \emph{The Riemann Zeta function and coin tossing}, {\sc Journal f{\"u}r die reine und angewandte Mathematik (Crelle's Journal)}, \textbf{254} (1972), pp. 100--109

  %\bibitem{Fujii} A. Fujii, \emph{Uniform Distribution of the Zeros of the Riemann Zeta Function and the Mean Value Theorems of diricchlet L-functions}, {\sc Proc. Japan Acad.}, \textbf{63}, Ser. A (1987), 370--373

  \bibitem{Fujii2} A. Fujii, \emph{On a Theorem of Landau. II}, {\sc Proc. Japan Acad.}, \textbf{66}, Ser. A (1990), 291--296

    \bibitem{Hlawka}  E. Hlawka, \emph{\"{U}ber die Gleichverteilung gewisser Folgen, welche mit den Nullstellen der Zetafunktion zusammenh\"{a}ngen} {\sc \"{O}sterreich. Akad. Wiss. Math.- Naturwiss. Kl. S.-B. II}, \textbf{184} (1975), no. 8–10, 459--471

\bibitem{Knopp} M. Knopp and S. Robins, \emph{Easy proofs of the Riemann's functional equation for $\zeta(s)$ and of Lipschitz summation}, {\sc Proc. AMS}, \textbf{129} (2001), pp. 1915--1922

\bibitem{Landau} E. Landau, \emph{\"{U}ber die Nullstellen der Zetafunktion}, {\sc Math. Ann.} \textbf{71} (1912), pp. 548--564

\bibitem{Montgomery} H. L. Montgomery, \emph{The Pair Correlation of Zeros of the Zeta Function}, In: H. G. Diamond, editor,
\emph{Analytic Number Theory, Proc. Symp. Pure Math.},  American Mathematical Society, Providence (1973), pp. 181--193

\bibitem{review} H. Niederriter, AMS MathScieNet review MR0453661
%\url{http://www.ams.org/mathscinet/search/publdoc.html?pg1=MR&s1=0453661&loc=fromreflist} (1978)

\bibitem{Odlyzko} A. Odlyzko, personal website, \url{http://www.dtc.umn.edu/~Odlyzko/zeta_tables/index.html} (2015)

\bibitem{Rademacher} H.A. Rademacher, \emph{Fourier Analysis in Number Theory, Symposium on Harmonic Analysis and Related Integral Transforms (Cornell Univ., Ithaca, N.Y., 1956)} in: Collected Papers of Hans Rademacher, Vol. II, pp. 434--458, Massachusetts Inst. Tech., Cambridge, Mass., 1974

\bibitem{sowa21} A. Sowa, \textit{The Dirichlet ring and unconditional bases in $L_2[0,2\pi]$}, {\sc Func. Anal. Appl.}, \textbf{47} (2013), 227-232%; Transl. from {Funkts. Analiz i Ego Prilozheniya, Vol. 47, No. 3, pp. 75-81, 2013

\bibitem{sowa22} A. Sowa, \emph{Factorizing matrices by Dirichlet multiplication}, {\sc Lin. Alg. Appl.}, \textbf{438} (2013), 2385-2393

\bibitem{Steuding} J. Steuding, \emph{One Hundred Years Uniform Distribution Modulo One and Recent Applications to Riemann's Zeta-Function}, Topics in mathematical analysis and applications, {\sc Springer Optim. Appl.} \textbf{94} (2014), pp. 659--698



\end{thebibliography}
\end{document}